\documentclass[11pt]{amsart}
\usepackage{amsmath,amsfonts, amssymb,amsthm}

\def\E{{\bf E}}

\def\Var{{\rm Var}}

\def\be{\begin{equation}}
\def\en{\end{equation}}

\def\bee{\begin{eqnarray*}}
\def\ene{\end{eqnarray*}}

\def\ep{\varepsilon}
\def\R{{\bf R}}

\title{Bounds on the deficit in the logarithmic Sobolev inequality}

\address{School of Mathematics, University of Minnesota, USA} 
\email{bobkov@math.umn.edu}

\address{Universit\'e Paris Est Marne la Vall\'ee - Laboratoire d'Analyse 
et de Math\'e\-matiques Appliqu\'ees (UMR CNRS 8050), 5 bd Descartes, 
77454 Marne la Vall\'ee Cedex 2, France}
\email{nathael.gozlan@univ-mlv.fr, paul-marie.samson@univ-mlv.fr}

\address{Universit\'e Paris Ouest Nanterre la D\'efense, MODAL'X, 
EA 3454, 200 avenue de la R\'epublique 92000 Nanterre, France}
\email{croberto@math.cnrs.fr}

\thanks{The first author was partially supported by NSF grant DMS-1106530. 
The other authors were partially supported by the ``Agence Nationale 
de la Recherche'' through the grants ANR 2011 BS01 007 01,
ANR-10-LABX-58 and ANR-11-LABX-0023-01.}

\vspace{2in}

\author{S. G. Bobkov, N. Gozlan, C. Roberto and P.-M. Samson}

\begin{document}
\bibliographystyle{plain}

\begin{abstract}

The deficit in the logarithmic Sobolev inequality for the Gaussian measure
is considered and estimated by means of transport and information-theoretic 
distances.
\end{abstract}
\maketitle

% -----------------------------  section 1 ---------------------------

\section{Introduction}
\setcounter{equation}{0}

Let $\gamma$ denote the standard Gaussian measure on the Euclidean space 
$\R^n$, thus with density
$$
\frac{d\gamma(x)}{dx} = \frac{1}{(2\pi)^{n/2}}\,e^{-|x|^2/2}
$$
with respect to the Lebesgue measure. (Here and in the sequel $|x|$ stands 
for the Euclidean norm of a vector $x \in \R^n$.) One of the basic results 
in the Gaussian Analysis is the celebrated logarithmic Sobolev inequality
\be
\int f\log f\,d\gamma - \int f\,d\gamma \, \log \int f\,d\gamma 
\leq \frac{1}{2}\, \int \frac{|\nabla f|^2}{f}\,d\gamma,
\en
holding true for all positive smooth functions $f$ on $\R^n$ with gradient 
$\nabla f$. In this explicit form it was obtained in the work of L. Gross [G],
initiating fruitful investigations around logarithmic Sobolev inequalities
and their applications in different fields. See e.g. a survey by M. Ledoux [L1]
and the books [L2,A] for a comprehensive account of such activities up to 
the end of 90's. One should mention that in an equivalent form -- as 
a relation between the Shannon entropy and the Fisher information, 
(1.1) goes back to the work by A. J. Stam [St], see [A, Chapter 10].

The inequality (1.1) is homogeneous in $f$, so the restriction 
$\int f\,d\gamma = 1$ does not lose generality. It is sharp in the sense 
that the equality is attained, namely for all $f(x) = e^{l(x)}$ with 
arbitrary affine functions $l$ on $\R^n$ (in which case the measures 
$\mu = f\gamma$ are still Gaussian). It is nevertheless of a certain interest 
to realize how large the difference between both sides of (1.1) is. 
This problem has many interesting aspects. For example, as was shown by 
E. Carlen in [C], which was perhaps a first address of the sharpness problem, 
for $f = |u|^2$ with a smooth complex-valued $u$ such that 
$\int |u|^2\,d\gamma = 1$, (1.1) may be strengthened to
$$
\int |u|^2\log |u|^2\,d\gamma + \int |W u|^2\log |W u|^2\,d\gamma \leq 
2 \int |\nabla u|^2\,d\gamma,
$$
where $W$ denotes the Wiener transform of $u$. That is, a certain non-trivial
functional may be added to the left-hand side of (1.1).

One may naturally wonder how to bound from below the deficit in (1.1), 
that is, the quantity
$$
\delta(f) = \frac{1}{2}\,\int \frac{|\nabla f|^2}{f}\,d\gamma -
\bigg[\int f\log f\,d\gamma - \int f\,d\gamma \, \log \int f\,d\gamma\bigg],
$$
in terms of more explicit, like distribution-dependent characteristics of $f$ 
showing its closeness to the extremal functions $e^{l}$ (when $\delta(f)$ 
is small). Recently, results of this type have been obtained by A. Cianchi, 
N. Fusco, F. Maggi and A. Pratelli [C-F-M-P] in their study of the closely 
related isoperimetric inequality for the Gaussian measure. The work by 
E. Mossel and J. Neeman [M-N] deals with dimension-free bounds for the 
deficit in one functional form of the Gaussian isoperimetric inequality 
appearing in [B]. See also the subsequent paper by R. Eldan [E] where
almost tight two-sided robustness bounds have been derived.
In [F-M-P1,Se] the authors deal with quantitative Brunn-Minkowski inequality (which is related to
the isoperimetric problem in Euclidean space), while bounds on the deficit in the Sobolev inequalities can be found in e.g.\ [F-M-P2,D-T] and in the Gagliardo-Nirenberg-Sobolev inequality in [C-F] (see also the references therein for more on the literature).

As for (1.1), one may also want to involve distance-like quantities between 
the measures $\mu = f\gamma$ and $\gamma$. This approach looks even more 
natural, when the logarithmic Sobolev inequality is treated as the relation 
between classical information-theoretic distances as 
\be
D(X|Z) \leq \frac{1}{2}\,I(X|Z).
\en

To clarify this inequality, let us recall standard notations and definitions.
If random vectors $X$ and $Z$ in $\R^n$ have distributions $\mu$ and $\nu$ 
with densities $p$ and $q$, and $\mu$ is absolutely continuous with respect
to $\nu$, the relative entropy of $\mu$ with respect to $\nu$ is defined by
$$
D(X|Z) = D(\mu|\nu) = \int p(x)\,\log\frac{p(x)}{q(x)}\,dx.
$$
%(and $D(X|Z) = \infty$ in all other cases). 
Moreover, if $p$ and $q$ are smooth, one defines the relative Fisher information
$$
I(X|Z) = I(\mu|\nu) = \int 
\Big|\frac{\nabla p(x)}{p(x)} - \frac{\nabla q(x)}{q(x)}\Big|^2\,p(x)\,dx.
$$
Both quantities are non-negative, and although non-symmetric in $(\mu,\nu)$, 
they may be viewed as strong distances of $\mu$ to $\nu$. This is already 
demonstrated by the well-known Pinsker inequality [P], connecting $D$ with 
the total variation norm:
$$
D(\mu|\nu) \geq \frac{1}{2}\,\|\mu - \nu\|_{\rm TV}^2.
$$

In the sequel, we mainly consider the particular case where $Z$ is standard 
normal, so that $\nu = \gamma$ in the above formulas. And in this case, 
as easy to see, for $d\mu = f\,d\gamma$ with $\int f\,d\gamma = 1$, 
the logarithmic Sobolev inequality (1.1) turns exactly into (1.2).

The aim of this note is to develop several lower bounds on the deficit in 
this inequality, $\frac{1}{2}\,I(X|Z) - D(X|Z)$, by involving also transport 
metrics such as the quadratic Kantorovich distance (see e.g.\ [V])
$$
W_2(X,Z) = W_2(\mu,\gamma) = 
\inf_\pi \Big(\int\!\!\!\int |x-z|^2\,d\pi(x,z)\Big)^{1/2}
$$
(where the infimum runs over all probability measures on $\R^n \times \R^n$ 
with marginals $\mu$ and $\gamma$). More generally, one may consider the 
optimal transport cost
$$
\mathcal{ T}(X,Z) = \mathcal{T}(\mu,\gamma) = 
\inf_\pi \int\!\!\!\int c(x-z)\,d\pi(x,z)
$$
for various ``cost" functions $c(x-z)$.

The metric $W_2$ is of weak type in the sense that it metrizes 
the weak topology in the space of probability measures on $\R^n$ 
(under proper moment constraints). It may be connected with the relative 
entropy by virtue of M. Talagrand's transport-entropy inequality
\be
W_2(X,Z)^2 \leq 2\,D(X|Z),
\en
cf.\ [T]. In view of (1.2), this also gives
an apriori weaker transport-Fisher information inequality
\be
W_2(X,Z) \leq \sqrt{I(X|Z)}.
\en

In formulations below, we use the non-negative convex function 
$$
\Delta(t)=t-\log(1+t), \qquad t > -1,
$$ 
and denote by $Z$ a random vector in $\R^n$ with the standard normal law.
%Let us now state a few results. 

\vskip5mm
{\bf Theorem 1.1.} {\it For any random vector $X$ in $\R^n$ with a smooth 
density, such that $I(X|Z)$ is finite,
\be
I(X|Z)-2D(X|Z) \, \geq \, n\Delta\Big(\frac{I(X)}{n} - 1\Big).
\en
Moreover,
\be
I(X|Z)-2D(X|Z) \, \geq \, \big(\sqrt{I(X|Z)}-W_2(X,Z)\big)^2 + 
n\Delta\left(\frac{W_2(X,Z)}{\sqrt{I(X|Z)}}\, 
\Big(\frac{I(X)}{n} - 1\Big)\!\!\right).
\en
}

\vskip5mm
As is common,
$$
I(X) = \int \frac{|\nabla p(x)|^2}{p(x)}\,dx
$$
stands for the usual (non-relative) Fisher information. Thus, (1.5)-(1.6) 
represent certain sharpenings of the logarithmic Sobolev inequality.
The lower bounds of the deficit in (1.5) and (1.6) are not simply comparable. 
However, in the next section, we recall that (1.5) is a self improvement 
of the logarithmic Sobolev inequality that obviously follows from (1.6).

An interesting feature of the bound (1.6) is that, by removing the last term 
in it, we arrive at the Gaussian case in the so-called HWI inequality due 
to F. Otto and C. Villani [O-V],
\be
D(X|Z)\leq W_2(X,Z)\sqrt{I(X|Z)} - \frac{1}{2}\,W_2^2(X,Z).
\en

As for (1.5), its main point is that, when $\E\,|X|^2 \leq n$, then
necessarily $I(X) \geq n$, and moreover, one can use the lower bound
$$
\frac{1}{n}\,I(X) - 1 = \frac{1}{n}\,I(X|Z) - \frac{1}{n}\,\E\,|X|^2 + 1 \geq 
\frac{1}{n}\,I(X|Z).
$$
Since $\Delta(t)$ is increasing for $t \geq 0$, (1.5) is then simplified to 
\be
I(X|Z)-2D(X|Z) \geq n\Delta\Big(\frac{1}{n}\,I(X|Z)\Big).
\en
In fact, this estimate is rather elementary in that it surprisingly follows 
from the logarithmic Sobolev inequality itself by virtue of rescaling 
(as will be explained later on). Here, let us only stress that the right-hand 
side of (1.8) can further be bounded from below. For example, by (1.2)-(1.3), 
we have
$$
I(X|Z)-2D(X|Z)  \, \geq \, n\Delta\Big(\frac{2}{n}\,D(X,Z)\Big)
\, \geq \, n\Delta\Big(\frac{1}{n}\,W_2^2(X,Z)\Big).
$$
But, $\frac{1}{n}\,W_2^2(X,Z) \leq \frac{1}{n}\,\E\,|X-Z|^2 \leq 4$, and
using $\Delta(t) \sim \frac{t^2}{2}$ for small $t$, the above yields a simpler 
bound.

\vskip5mm
{\bf Corollary 1.2.} {\it For any random vector $X$ in $\R^n$ with a smooth 
density and such that $\E\,|X|^2 \leq n$, we have
\be
I(X|Z)-2D(X|Z) \, \geq \, \frac{c}{n}\,W_2^4(X,Z),
\en
up to an absolute constant $c>0$.
}

\vskip5mm
{\bf Remark.} Dimensional refinements of the HWI inequality (1.7) similar to
(1.6) were recently considered by several authors. For instance, F-Y. Wang 
obtained in [W] some HWI type inequalities involving the dimension and the 
quadratic Kantorovich distance under the assumption that the reference measure 
enjoys some curvature dimension condition $\mathrm{CD}(-K,N)$ with $K\geq0$ 
and $N\geq0$ (see [B-E] for the definition). See also the recent paper [E-K-S] 
for dimensional variants of the HWI inequality in an abstract metric space 
framework. The standard Gaussian measure does not enter directly the framework 
of [W] (or [E-K-S]), but we believe that it might be possible to use similar 
semigroup arguments to derive (1.6). In the same spirit, D. Bakry, F. Bolley 
and I. Gentil [B-B-G]  used semigroup techniques to prove a dimensional 
reinforcement of Talagrand's transport-entropy inequality. 

\vskip2mm
Returning to (1.9), we note that, after a certain recentering of $X$, one 
may give some refinement over this bound, especially when $D(X|Z)$ is small. 
Given a random vector $X$ in $\R^n$ with finite absolute moment, define 
the recentered random vector $\bar X = (\bar X_1,\dots,\bar X_n)$ by putting 
$\bar X_1 = X_1 - \E X_1$ and
$$
\bar X_k = X_k - \E\, (X_k | X_1,\dots,X_{k-1}), \qquad k \geq 2,
$$
where we use standard notations for the conditional expectations.

\vskip5mm
{\bf Theorem 1.3.} {\it For any random vector $X$ in $\R^n$ with a smooth 
density, such that $I(X|Z)$ is finite, the deficit in $(1.2)$ satisfies
\be
\frac{1}{2}\,I(X|Z) - D(X|Z) \geq c\, 
\frac{\mathcal{T}^2(\bar X,Z)}{D(\bar X|Z)}.
\en
Here the optimal transport cost $\mathcal{T}$ corresponds to the cost function 
$\Delta(|x-z|)$, $c$ is a positive absolute constant and one uses the 
convention $0/0=0$ in the right hand side.
}

\vskip5mm
In particular, in dimension one, if a random vector $X$ has mean zero, we get
that
\be
\frac{1}{2}\,I(X|Z) - D(X|Z) \geq c\, \frac{\mathcal{T}^2(X,Z)}{D(X|Z)}.
\en

The bound (1.10) allows one to recognize the cases of equality in 
(1.2) -- this is only possible when the random vector $X$ is a translation of 
the standard random vector $Z$ (an observation of E. Carlen [C] who used 
a different proof). The argument is sketched in Appendix C.
%Indeed, if $X$ is not normal, then necessarily $\mathcal{T}(\bar X,Z) > 0$.

It is worthwhile noting that the transport cost $\mathcal{T}$ of Theorem 1.3
already appeared in the literature, cf. e.g. [B-G-L] or [B-K]. In particular, 
it was shown in [B-G-L] that this transport cost can be used to give 
an alternative representation of the Poincar\'e inequality. In fact, 
it may be connected with the classical Kantorovich transport distance $W_1$ 
based on the cost function $c(x,z) = |x-z|$. More precisely, due to the 
convexity of $\Delta$, there are simple bounds
$$
W_1(X,Z) \, \geq \,
\mathcal{T}(X,Z) \, \geq \, \Delta(W_1(X,Z)) \sim \min\{W_1(X,Z),W_1^2(X,Z)\}.
$$
Hence, if $D(\bar X|Z) \leq 1$, then according to (1.3),
$W_1^2(X,Z) \leq W_2^2(X,Z) \leq 2$, and (1.10) is simplified to
\be
\frac{1}{2}\,I(X|Z) - D(X|Z) \geq c'\, \frac{W_1^4(\bar X,Z)}{D(\bar X|Z)},
\en
for some other absolute constant $c'.$

In connection with such bounds, let us mention a recent preprint by
E. Indrei and D. Marcon [I-M], which we learned about while the current 
work was in progress. For a $C^2$-smooth function $V$ on $\R^n$, let us 
denote by $V''(x)$ the matrix of second partial derivatives of $V$ at 
the point $x$. We use comparison of symmetric matrices in the usual matrix 
sense and denote by ${\rm I_n}$ the identity $n \times n$ matrix.

It is proved in [I-M] (Theorem 1.1 and Corollary 1.2) that, if a random 
vector $X$ on $\R^n$ has a smooth density $p = e^{-V}$ satisfying 
$\ep\,{\rm I_n} \leq V'' \leq M\,{\rm I_n}$ ($0<\ep<M$), then
\be
\frac{1}{2}\,I(X|Z)-D(X|Z) \, \geq \, c\,W_2^2(X-\E X,Z)
\en
with some constants $c = c(\ep,M)$. In certain cases it is somewhat stronger 
than (1.11). We will show that a slight adaptation of our proof 
of (1.11) leads to a bound similar to (1.13).

\vskip5mm
{\bf Theorem 1.4.} {\it Let $X$ be a random vector in $\R^n$ with a smooth 
density $p = e^{-V}$ with respect to Lebesgue measure such that
$V'' \geq \ep\,{\rm I_n}$, for some $\ep>0$. Then, the deficit in $(1.2)$ 
satisfies
\be
\frac{1}{2}\,I(X|Z) - D(X|Z) \geq c\, \min(1,\ep)\, W_2^2(\bar{X},Z),
\en
for some absolute constant $c$.
}

\vskip5mm
Note that Theorem 1.4 holds under less restrictive assumptions on $p$ 
than the result from [I-M]. In particular, in dimension $1$, we see that 
the constant $c$ in (1.13) can be taken independent on $M$.
In higher dimensions however, it is not clear how to compare $W_2(\bar{X},Z)$ 
and $W_2(X-\E X, Z)$ in general. One favorable case is, for instance, when 
the distribution of $X$ is unconditional (i.e., when its density $p$ 
satisfies $p(x)=p(\ep_1x_1,\dots,\ep_nx_n)$, for all $x \in \R^n$ and all 
$\ep_i = \pm 1$). In this case, $\E X = 0$ and $\bar{X} = X$, and 
thus (1.14) reduces to (1.13) with a constant $c$ independent on $M$.

Let us mention that in Theorem 1.3 of [I-M], the assumption 
$V''\leq M\,{\rm I_n}$ can be relaxed into an integrability condition of 
the form $\int \|V''\|^r\,dx \leq M$, for some $r>1$, but only at the expense 
of a constant $c$ depending on the dimension $n$ and of an exponent greater 
than $2$ in the right-hand side of (1.13).

Finally, let us conclude this introduction by showing optimality of the 
bounds (1.11), (1.12), (1.14) for mean zero Gaussian random vectors with 
variance close to $1$. An easy calculation shows that, if $Z$ is a standard 
Gaussian random vector in $\R^n$, then for any $\sigma>0$,
$$
D(\sigma Z|Z)=\frac{n}{2}\left((\sigma^2-1)-2\log \sigma \right), \qquad 
I(\sigma Z|Z)= n\sigma^2\Big(\frac{1}{\sigma^2}-1\Big)^2,
$$
so that
$$
\frac{1}{2}\,I(X|Z)-D(X|Z) \, = \, \frac{n}{2}\,
\Big(\frac{1}{\sigma^2} - 1 + 2\log \sigma\Big) \, \sim \, n(\sigma-1)^2, 
\quad \ \ \text{as} \ \ \sigma \to 1.
$$
On the other hand,
$$\quad W_2^2(\sigma Z,Z) = n (\sigma-1)^2, \qquad 
W_1(\sigma Z,Z)=|\sigma-1|\, \E\, |Z| \, \simeq \, |\sigma-1|\sqrt{n},
$$
and thus the three quantities 
$W_2^2(\sigma Z,Z)$, $\mathcal{T}^2(\sigma Z,Z)/D(\sigma Z |Z)$ and 
$W_1^4(\sigma Z,Z)/D(\sigma Z|Z)$ are all of the same order $n(\sigma-1)^2$, 
when $\sigma$ goes to $1$.

\vskip5mm
The paper is organized in the following way. In Section 2 we recall
Stam's formulation of the logarithmic Sobolev inequality in the form of
an ``isoperimetric inequality for entropies'' and discuss the involved
improved variants of (1.1). Theorem 1.1 is proved in Section 3.
In Section 4 we consider sharpened transport-entropy inequalities 
in dimension one, which are used to derive bounds on the deficit like 
those in (1.11)-(1.14). For general dimensions Theorems 1.3 and 1.4 are 
proved in Section 5. For the reader's convenience and so as to get 
a more self-contained exposition, we move to Appendices several 
known results and arguments.

% -----------------------------  section 2 ---------------------------

\vskip10mm
\section{Self-improvement of the logarithmic Sobolev inequality}
\setcounter{equation}{0}

To start with, let us return to the history and remind the reader Stam's
information-theoretic formulation of the logarithmic Sobolev inequality.
As a base for the derivation, one may take (1.2) and rewrite it in terms of 
the Fisher information $I(X)$ and the (Shannon) entropy
$$
h(X) = - \int p(x)\,\log p(x) \,dx,
$$
where $X$ is a random vector in $\R^n$ with density $p$. Here the integral 
is well-defined, as long as $X$ has finite second moment. Introduce also 
the entropy power
$$
N(X) = \exp\{2h(X)/n\},
$$
which is a homogeneous functional of order 2. The basic connections between 
the relative and non-relative information quantities are given by
$$
D(X|Z) = h(Z) - h(X), \qquad I(X|Z) = I(X) - I(Z),
$$
where $Z$ has a normal distribution, and provided that $\E\,|X|^2 = \E\,|Z|^2$.

More generally, assuming that $Z$ is standard normal and 
$\E\,|X|^2 < \infty$, the first above equality should be replaced with
$$
D(X|Z) = - h(X) + \E\, \Big(\frac{n}{2}\,\log(2\pi) + \frac{|X|^2}{2}\Big),
$$
while, as was mentioned before, under mild regularity assumptions on $p$, 
$$
I(X|Z) = I(X) + \E\,|X|^2 -2n.
$$
Inserting these expressions into the inequality (1.2),
the second moment is cancelled, and (1.2) becomes
$$
I(X) + 2h(X) \geq 2n + n\log(2\pi).
$$
However, this inequality is not homogeneous in $X$. So, one may apply
it to $\lambda X$ in place of $X$ with arbitrary $\lambda>0$ and then 
optimize. The function
$$
v(\lambda) = I(\lambda X) + 2 h(\lambda X) = 
\frac{I(X)}{\lambda^2} + n\log \lambda^2 + 2 h(X)
$$
is minimized for $\lambda^2 = I(X)/n$, and at this point the inequality
becomes:

\vskip5mm
{\bf Theorem 2.1} ([St]). {\it If a random vector $X$ in $\R^n$ has a smooth
density and finite second moment, then
\be
I(X) \, \frac{N(X)}{2\pi e} \, \geq \, n.
\en
}

\vskip5mm
This relation was first obtained by Stam and is sometimes referred 
to as the isoperimetric inequality for entropies, cf.\ e.g.\ [D-C-T]. 
Stam's original argument is based on the general entropy power inequality
\be
N(X+Y) \geq N(X) + N(Y),
\en
which holds for all independent random vectors $X$ and $Y$ in $\R^n$ with 
finite second moments (so that the involved entropies do exist, cf.\ also 
[Bl], [Li]). Then, (2.1) can be obtained by taking $Y=\sqrt{t}\,Z$ with $Z$ 
having a standard normal law and combining (2.2) with the de Bruijn identity
\be
\frac{d}{dt}\, h(X+\sqrt{t}\,Z) = \frac{1}{2}\,I(X+\sqrt{t}\,Z) \qquad (t>0).
\en

Note that in the derivation $(1.2) \Rightarrow (2.1)$ the argument may easily 
be reversed, so these inequalities are in fact equivalent (as noticed by 
E. Carlen [C]). On the other hand, the isoperimetric inequality for entropies 
can be viewed as a certain sharpening of (1.1)-(1.2). Indeed, let us rewrite 
(2.1) explicitly as
\be
\int p(x)\,\log p(x) \,dx \, \leq \, \frac{n}{2}\,
\log\Big(\frac{1}{2\pi e\, n}\,\int \frac{|\nabla p(x)|^2}{p(x)}\,dx\Big).
\en
It is also called an optimal Euclidean logarithmic Sobolev inequality;
cf. [B-L] for a detail discussion including deep connections with
dimensional lower estimates on heat kernel measures. In terms of the density 
$f(x) = \sqrt{2\pi}e^{x^2/2}p(x)$ of $X$ with respect to $\gamma$ we have
$$
\int p(x)\,\log p(x) \,dx = \frac{n}{2}\,\log\frac{1}{2\pi} - 
\frac{1}{2}\,\int |x|^2 f(x)\,d\gamma(x) + \int f\,\log f\,d\gamma,
$$
while
$$
\int \frac{|\nabla p(x)|^2}{p(x)}\,dx = 
\int \frac{|\nabla f(x)|^2}{f(x)}\,d\gamma(x) - 
\int |x|^2 f(x)\,d\gamma(x) + 2n.
$$
Inserting these two equalities in (2.4), we arrive at the following 
reformulation of Theorem 2.1.

\vskip5mm
{\bf Corollary 2.2.} {\it For any positive smooth function $f$ on $\R^n$ 
such that $\int f\,d\gamma = 1$, putting 
$b = \frac{1}{n}\int |x|^2 f(x)\,d\gamma(x)$, we have
\be
\int f\,\log f\,d\gamma \, \leq \, \frac{n}{2}\,
\log\Big(\frac{1}{n}\int \frac{|\nabla f|^2}{f}\,d\gamma + (2-b)\Big) + \frac{n}{2}\,(b-1),
\en
which is exactly (1.5). In particular, if $b \leq 1$,
\be
\int f\,\log f\,d\gamma \, \leq \, \frac{n}{2}\,
\log\Big(\frac{1}{n}\int \frac{|\nabla f|^2}{f}\,d\gamma + 1\Big).
\en
}

\vskip5mm
An application of $\log t \leq t-1$ on the right-hand side of (2.5) returns 
us to the original logarithmic Sobolev inequality (1.1). It is in this sense that
Inequality (2.5) is stronger, although it was derived from
(1.1). The point of self-improvement is that the $\log$-value of 
$$
I = \int \frac{|\nabla f|^2}{f}\,d\gamma
$$ 
may be much smaller than the integral itself. This can be used, for example, 
in bounding the deficit $\delta(f)$ in (1.1). Indeed, when $b \leq 1$, (2.6) 
yields
$$
2\delta(f) \geq I - n\,\log\Big(\frac{1}{n}\,I + 1\Big).
$$
That is, using again the function $\Delta(t) = t - \log(t+1)$, we have
$$
2\delta(f) \geq 
n\,\Delta\Big(\frac{1}{n}\int \frac{|\nabla f|^2}{f}\,d\gamma\Big).
$$
But this is exactly the information-theoretic bound (1.8), mentioned 
in Section 1 as a direct consequence of (1.5).

As the function $\Delta$ naturally appears in many related inequalities,
let us collect together a few elementary bounds that will be needed in 
the sequel.

\vskip5mm
{\bf Lemma 2.3.} {\it We have:

\vskip2mm
$a)$ \ $\Delta(ct) \geq \min(c,c^2)\,\Delta(t)$, whenever $c,t \geq 0$;

\vskip1mm
$b)$ \ $\Delta(t) \geq \frac{1}{2}\,t^2$, for all $-1 < t \leq 0$;

\vskip1mm
$c)$ \ $\Delta(t) \geq \frac{\Delta(a)}{a^2}\,t^2$, for all $0 \leq t \leq a$ 
$(a>0)$;

\vskip1mm
$d)$ \ $(1 - \log 2)\,\min\{t,t^2\} \leq \Delta(t) \leq t$, for all $t \geq 0$. 

\vskip2mm
Moreover, for any random variable $\xi \geq 0$,
$$
(1 - \log 2)\,\min\{\E \xi, (\E \xi)^2\} \leq \E \Delta(\xi) \leq \E \xi.
$$

}

\vskip2mm
{\bf Proof.} $a)$
In case $0 \leq c \leq 1$, the required inequality follows from the 
representation 
$$
\Delta(ct) = \int_0^{ct} \Delta'(s)\,ds = \int_0^{ct} \frac{s}{1+s}\,ds = 
c^2\int_0^{t} \frac{u}{1+cu}\,du.
$$
In case $c \geq 1$, it becomes $\log(1 + ct) \leq c\log(1+t)$, which is
obvious.

$b)$ This bound immediately follows from the Taylor expansion for 
the function $-\log(1-s)$.

$c)$ It is easy to check that the function $\Delta(\sqrt{x})$ is concave
in $x \geq 0$. Hence, the optimal value of the constant $c$ in
$\Delta(t) \geq ct^2$ on the interval $[0,a]$ corresponds to the endpoint 
$t = a$.

$d)$ For $t \geq 1$, the first inequality becomes 
$ct \leq t - \log(1+t)$, where $c = 1 - \log 2$. Both sides are equal 
at $t=1$, and we have inequality for the derivatives at this point. 
Hence, it holds for all $t \geq 1$. For the interval $0 \leq t \leq 1$,
the inequality $\Delta(t) \geq ct^2$ is given in $c)$.

Finally, an application of Jensen's inequality with the convex function 
$\Delta$ together with $\Delta(\xi) \leq \xi$
leads to the last bounds of the lemma.
\qed

% -----------------------------  section 3 ---------------------------

\vskip10mm
\section{HWI inequality and its sharpening}
\setcounter{equation}{0}

We now turn to the remarkable HWI inequality of F. Otto and C. Villani and 
state it in full generality. Assume that the probability measure $\nu$ 
on $\R^n$ has density
$$
\frac{d\nu(x)}{dx} = e^{-V(x)}
$$
with a twice continuously differentiable $V:\R^n \rightarrow \R$. 

\vskip5mm
{\bf Theorem 3.1} ([O-V]). {\it \, Assume that $V''(x) \geq \kappa\,{\rm I_n}$ 
for all $x \in \R^n$ with some $\kappa \in \R$. Then, for any probability 
measure $\mu$ on $\R^n$ with finite second moment,
\be
D(\mu|\nu) \leq 
W_2(\mu|\nu) \sqrt{I(\mu|\nu)} - \frac{\kappa}{2}\,W_2^2(\mu,\nu).
\en
}

\vskip5mm
This inequality connects together all three important distances: the relative 
entropy (which sometimes is denoted by $H$), the relative Fisher information 
$I$, and the quadratic transport distance $W_2$. It may equivalently be 
written as
\be
D(\mu|\nu) \leq 
\frac{1}{2\ep}\, I(\mu|\nu) + \frac{\ep - \kappa}{2}\,W_2^2(\mu,\nu)
\en
with an arbitrary $\ep > 0$. Taking here $\ep = \kappa$, one gets
\be
D(\mu|\nu) \leq \frac{1}{2\kappa}\, I(\mu|\nu).
\en
If $\nu = \gamma$, we arrive in (3.3) at the logarithmic Sobolev inequality 
(1.1) for the Gaussian measure, and thus the HWI inequality represents its 
certain refinement. In particular, (3.1) may potentially be used in the study 
of the deficit in (1.1), as is pointed in Theorem 1.1. 

In the proof of the latter, we will use two results. The following lemma, 
reversing the transport-entropy inequality, may be found in the survey by 
Raginsky and Sason [R-S], Lemma 15. It is due to Y. Wu [Wu] who used it 
to prove a weak version of the Gaussian HWI inequality (without 
the curvature term $-\frac{1}{2}W_2^2(X,Z)$ appearing in (1.7)). 
The  proof of Lemma 3.2 is reproduced in Appendix A. 

For a random vector $X$ in $\R^n$ with finite second moment, put
$$
X_t = X + \sqrt{t}\,Z \qquad (t \geq 0),
$$
where $Z$ is a standard normal random vector in $\R^n$, independent of $X$.

\vskip5mm
{\bf Lemma 3.2.} {([Wu]) \it Given  random vectors $X$ and $Y$ in $\R^n$ 
with finite second moments, for all $t > 0$, 
$$
D(X_t |Y_t)\leq \frac{1}{2t}\, W_2^2(X,Y).
$$
}

\vskip5mm
We will also need a convexity property of the Fisher information in the form 
of the Fisher information inequality. As a full analog of the entropy power 
inequality (2.2), it was apparently first mentioned by Stam [St].

\vskip5mm
{\bf Lemma 3.3.} {\it Given independent random vectors $X$ and $Y$ in $\R^n$
with smooth densities,
\be
\frac{1}{I(X+Y)} \geq \frac{1}{I(X)} + \frac{1}{I(Y)}.
\en
}

\vskip5mm
{\bf Proof of Theorem 1.1.} Let $Z$ be standard normal, and let the 
distribution of $X$ not be a translation of $\gamma$ (in which case
both sides of (1.5) and of (1.6) are vanishing).

We recall that, if $Y$ is a normal random vector with mean zero and 
covariance matrix $\sigma^2\, {\rm I}_n$, then 
$$
D(X|Y) = h(Y)-h(X)+ \frac{1}{2\sigma^2}\left(\E\,|X|^2 - \E\,|Y|^2\right).
$$
In particular,
$$
D(X|Z) = h(Z)-h(X)+ \frac{1}{2}\left(\E\,|X|^2 - \E\,|Z|^2\right),
$$
where $\E\,|Z|^2 = n$. Using de-Bruijn's identity (2.3), 
$\frac{d}{dt}\, h(X_t) = \frac{1}{2}\,I(X_t)$, we therefore obtain that, 
for all $t>0$,
\bee
D(X_t|Z_t)
 & = &
h(Z_t) - h(X_t) + \frac{1}{2(1+t)} \left(\E\,|X_t|^2 - \E\, |Z_t|^2\right) \\
 & = &
h(Z_t) - h(X_t) + \frac{1}{2(1+t)} \left(\E\,|X|^2 - \E\,|Z|^2\right) \\
 & = &
(h(Z) - h(X)) +  \frac{1}{2} \int_0^t (I(Z_\tau) - I(X_\tau))\,d\tau +
\frac{1}{2(1+t)}\left(\E\,|X|^2 - \E\,|Z|^2\right) \\
 & = &
D(X|Z) + \frac{1}{2} \int_0^t (I(Z_\tau) - I(X_\tau))\,d\tau -
\frac{t}{2(1+t)}\left(\E\,|X|^2 - \E\,|Z|^2\right).
\ene
Equivalently,
\be
D(X|Z) \, = \, D(X_t|Z_t) + 
\frac{1}{2} \int_0^t (I(X_\tau) - I(Z_\tau))\,d\tau +
\frac{t}{2(1+t)}\left(\E\,|X|^2 - \E\,|Z|^2\right).
\en
In order to estimate from above the last integral, we apply Lemma 3.3 to 
the couple $(X,\sqrt{\tau}\,Z)$, which gives
$$
I(X_\tau) \leq \frac{1}{\frac{1}{I(X)}+\frac{1}{I(\sqrt{\tau}\, Z)}} =
\frac{nI(X)}{n + \tau I(X)}.
$$
Inserting also $I(Z_\tau) = \frac{n}{1+\tau}$, we get
\bee
\int_0^t (I(X_\tau) - I(Z_\tau))\,d\tau 
 & \leq &
\int_0^t \left(\frac{nI(X)}{n + \tau I(X)} - \frac{n}{1+\tau}\right) \,d\tau \\ 
 & = &
\frac{n}{2}\log\frac{n + tI(X)}{n(1+t)}.
\ene
Thus, from (3.5),
$$
D(X|Z) \, \leq \, D(X_t|Z_t) + \frac{n}{2}\,\log\frac{n + tI(X)}{n(1 + t)} +
\frac{t}{2(1+t)}\left(\E\,|X|^2 - n\right).
$$
Furthermore, an application of Lemma 3.2 together with the identity 
$$
\E\,|X|^2 - n = I(X|Z) - I(X) + n
$$
yields
\be
D(X|Z) \, \leq \, \frac{1}{2t}\, W_2^2(X,Z) + 
\frac{n}{2}\,\log\frac{n + tI(X)}{n(1 + t)} +
\frac{t}{2(1+t)}\left(I(X|Z) - I(X) + n\right).
\en
As $t$ goes to infinity in (3.6), we get in the limit
$$
D(X|Z) \leq 
\frac{1}{2}\,I(X|Z) - \frac{n}{2}\,\Delta\Big(\frac{I(X)}{n} - 1\Big),
$$
which is exactly the required inequality (1.5) of Theorem 1.1.

As for (1.6), let us restate (3.6) as the property that the deficit 
$I(X|Z) - 2D(X|Z)$ is bounded from below by
\be
I(X|Z) - \frac{1}{t}\, W_2^2(X,Z) - 
n\,\log\frac{n + tI(X)}{n(1 + t)} - 
\frac{t}{1+t}\left(I(X|Z) - I(X) + n\right).
\en
Assuming that $X$ is not normal, we end the proof by choosing the value 
\be
t = \frac{W_2(X,Z)}{\sqrt{I(X|Z)} - W_2(X,Z)},
\en
which is well-defined and positive.
Indeed, by the assumption that $I(X|Z)$ is finite, $W_2(X,Z)$ 
is finite as well (according to the inequality (1.4), for example). 
Moreover, the case where $\sqrt{I(X|Z)} = W_2(X,Z)$ is impossible, since 
then $2D(X|Z) = I(X|Z)$. But the latter is only possible, when the 
distribution of $X$ represents a translation of $\gamma$, by the result of 
E. Carlen on the equality cases in (1.1) (cf. also Appendix C). 

Putting for short $W = W_2(X,Z)$, $I = I(X|Z)$, $I_0 = I(X)$, we finally
note that the expression (3.7) with the value of $t$ specified in 
(3.8) turns into
$$
\hskip-30mm
I - W(\sqrt{I} - W) - n\,\log\frac{1 + \frac{W}{\sqrt{I} - W}\, 
\frac{I_0}{n}}{\frac{\sqrt{I}}{\sqrt{I} - W}} - 
\frac{W}{\sqrt{I}}\left(I - I_0 + n\right)
$$
$$
\hskip30mm
= \ (\sqrt{I} - W)^2 - n\,
\log\Big(1 + \frac{W}{\sqrt{I}}\,\Big(\frac{I_0}{n} - 1\Big)\Big) + 
n\frac{W}{\sqrt{I}}\,\Big(\frac{I_0}{n} - 1\Big)
$$
$$
\hskip-8mm
= \ (\sqrt{I} - W)^2 + 
n\Delta\Big(\frac{W}{\sqrt{I}}\,\Big(\frac{I_0}{n} - 1\Big)\Big).
$$
\qed

% -----------------------------  section 4 ---------------------------

\vskip10mm
\section{Sharpened transport-entropy inequalities on the line}
\setcounter{equation}{0}

Nowadays, Talagrand's transport-entropy inequality (1.2),
\be
\frac{1}{2}\,W_2^2(\mu,\gamma) \leq D(\mu|\gamma),
\en
has many proofs (cf. e.g. [B-G]). In the one dimensional case it admits 
the following refinement, which is due to F. Barthe and A. Kolesnikov.

\vskip5mm
{\bf Theorem 4.1} ([B-K]). {\it For any probability measure $\mu$ on the 
real line with finite second moment, having the mean or median at the origin,
\be
\frac{1}{2}\,W_2^2(\mu,\gamma) + \frac{1}{4}\,\mathcal{T}'(\mu,\gamma) \leq D(\mu|\gamma),
\en
where the optimal transport cost $\mathcal{T}'$ is based on the cost function 
$c'(x-z) = \Delta\big(\frac{|x-z|}{\sqrt{2\pi}}\big)$.
}

\vskip5mm
It is also shown in [B-K] that the constant $\frac{1}{4}$ may be replaced 
with 1 under the median assumption. Anyhow, the deficit in (4.1) can be 
bounded in terms of the transport distance $\mathcal{T}$ which represents 
a slight weakening of $W_2$ (since the function $\Delta(t) = t - \log(t+1)$ 
is almost quadratic near zero).

In \cite{B-K}, the reinforced transport inequality above was only stated for probability measures with median at $0$, but the argument can be easily 
adapted to the mean zero case. For the sake of completeness, the proof of 
Theorem 4.1 is recalled in Appendix B. In order to work with the usual cost 
function $c(x-z) = \Delta(|x-z|)$, the inequality (4.2) will be modified to
\be
\frac{1}{2}\,W_2^2(\mu,\gamma) + \frac{1}{8\pi}\,\mathcal{T}(\mu,\gamma) 
\leq D(\mu|\gamma)
\en
under the assumption that $\mu$ has mean zero. (Here we use the elementary 
inequality $\Delta(ct)\geq c^2 \Delta(t)$, for $0\leq c\leq1$, $t\geq 0$, 
cf. Lemma 2.3.)

As a natural complement to Theorem 4.1, it will be also shown in Appendix B 
that, under an additional log-concavity assumption on $\mu$, the transport cost 
$\mathcal{T}$ in the inequalities (4.2)-(4.3) may be replaced with $W_2^2$.
That is, the constant $\frac{1}{2}$ in (4.1) may be increased.

\vskip5mm
{\bf Theorem 4.2.} {\it Suppose that the probability measure $\mu$ on the real 
line has a twice continuously differentiable density 
$\frac{d\mu}{dx}(x) = e^{-v(x)}$ such that, for a given $\ep>0$,
\be
v''(x) \geq \ep, \qquad x \in \R.
\en
If $\mu$ has mean at the origin, then with some absolute constant $c>0$ we have
\be
\Big(\frac{1}{2} + c\min\{1,\sqrt{\ep}\,\}\Big)\,W_2^2(\mu,\gamma) \leq 
D(\mu|\gamma).
\en
}

\vskip5mm
Here, one may take $c = 1-\log2$.

Let us now explain how these refinements can be used in the problem of bounding
the deficit in the one dimensional logarithmic Sobolev inequality.
Returning to (4.3), we are going to combine this bound with the HWI 
inequality (3.1). Putting 
$$
W = W_2(\mu,\gamma), \quad 
D = D(\mu|\gamma), \quad
I = I(\mu|\gamma), 
$$
we rewrite (3.1) as
$$
I - 2D \geq (\sqrt{I} - W)^2.
$$
On the other hand, applying the logarithmic Sobolev inequality $I \geq 2D$, 
(4.3) yields
$
I \geq W^2 + \frac{1}{4\pi}\,\mathcal{T},
$
where $\mathcal{T} = \mathcal{T}(\mu,\gamma)$. Hence,
$$
I - 2D \geq \bigg(\sqrt{W^2 + \frac{1}{4\pi}\,\mathcal{T}} - W\bigg)^2 =
W^2\, \bigg(\sqrt{1 + \frac{\mathcal{T}}{4\pi\,W^2}} - 1\bigg)^2.
$$
Here, by the very definition of the transport distance, one has
$\mathcal{T} \leq W^2$, so 
$\ep = \frac{\mathcal{T}}{4\pi\,W^2} \leq \frac{1}{4\pi}$. This implies that
$\sqrt{1 + \ep} - 1 \geq c\ep$ with 
$c = 4\pi\,\big(\sqrt{1 + \frac{1}{4\pi}} - 1\big)$.
Thus, up to a positive numerical constant,
\be
D + c\,\frac{\mathcal{T}^2}{W^2} \leq \frac{1}{2}\,I.
\en

In order to get a more flexible formulation, denote by $\mu_t$ the shift 
of the measure $\mu$,
$$
\mu_t(A) = \mu(A - t), \qquad A \subset \R \ \ ({\rm Borel}),
$$
which is the distribution of the random variable $X+t$ (with fixed $t \in \R$),
when $X$ has the distribution $\mu$. As easy to verify,
\bee
D(\mu_t|\gamma) & = & D(\mu|\gamma) + \frac{t^2}{2} + t\,\E X, \\
\frac{1}{2}\,I(\mu_t|\gamma)
 & = & 
\frac{1}{2}\,I(\mu|\gamma) + \frac{t^2}{2} + t\,\E X.
\ene
Hence, the deficit 
$$
\delta(\mu) = \frac{1}{2}\,I(\mu|\gamma) - D(\mu|\gamma)
$$
in the logarithmic Sobolev inequality (1.2) is translation invariant:
$\delta(\mu_t) = \delta(\mu)$. Applying (4.6) to $\mu_t$ with 
$t = -\int x\,d\mu(x)$, so that $\mu_t$ would have mean zero, therefore yields:

\vskip5mm
{\bf Corollary 4.3.} {\it For any non-Gaussian probability measure $\mu$ 
on the real line with finite second moment, up to an absolute constant $c>0$,
\be
D(\mu|\gamma) + c\,\frac{\mathcal{T}^2(\mu_{-t},\gamma)}{W_2^2(\mu_{-t},\gamma)} 
\leq
\frac{1}{2}\,I(\mu|\gamma),
\en
where the optimal transport cost $\mathcal{T}$ is based on the cost function 
$\Delta(|x-z|)$, and where $t= \int x \, d\mu(x)$. In particular,
\be
D(\mu|\gamma) + \frac{c}{2}\,\frac{\mathcal{T}^2(\mu_{-t},\gamma)}{D(\mu_{-t}|\gamma)} 
\leq
\frac{1}{2}\,I(\mu|\gamma).
\en
}

\vskip5mm
Here the second inequality follows from the first one by using 
$W_2^2 \leq 2D$. It will be used in the next section to perform tensorisation
for a multidimensional extension. Note that (4.8) may be derived directly
from (4.3) with similar arguments. Indeed, one can write
\bee
I - 2D
 & \geq &
(\sqrt{I} - W)^2 \ \geq \ (\sqrt{2D} - W)^2 \\
 & = &
\frac{(2D - W^2)^2}{(\sqrt{2D} + W)^2} \ \geq \
\frac{(2D - W^2)^2}{(2\sqrt{2D})^2} \ \geq \ 
\frac{\mathcal{T}^2}{128\,\pi^2 D^2},
\ene
thus proving (4.8) with constant 
$c = 1/(128\,\pi^2)$.

Let us now turn to Theorem 4.2 with its additional hypothesis (4.4).
Note that the property $v'' \geq 0$ describes the so-called log-concave 
probability distributions on the real line (with $C^2$-smooth densities), 
so (4.4) represents its certain quantitative strengthening. It is also 
equivalent to the property that $X$ has a log-concave density with respect 
to the Gaussian measure with mean zero and variance $\ep$.

Arguing as before, from (4.5) we have
$$
I - 2D \geq W^2\, \Big(\sqrt{1 + c\,\min\{1,\sqrt{\ep}\}} - 1\Big)^2.
$$
Hence, we obtain:

\vskip5mm
{\bf Corollary 4.4.} {\it Let $\mu$ be a probability measure on the real line 
with mean zero, and satisfying $(4.4)$ with some $\ep>0$. Then, 
up to an absolute constant $c>0$,
\be
D(\mu|\gamma) + c\,\min\{1,\ep\}\,W_2^2(\mu,\gamma) \leq 
\frac{1}{2}\,I(\mu|\gamma),
\en
}

% -----------------------------  section 5 ---------------------------

\vskip10mm
\section{Proof of Theorems~1.3 and 1.4}
\setcounter{equation}{0}

As the next step, it is natural to try to tensorize the inequality (4.8) 
so that to extend it to the multidimensional case. 

If $x = (x_1,\dots,x_n) \in \R^n$, denote by $x_{1:i}$ the subvector
$(x_1,\dots,x_i)$, $i = 1,\dots,n$. Given a probability measure $\mu$ on 
$\R^n$, denote by $\mu_1$ its projection to the first coordinate, i.e., 
$\mu_1(A) = \mu(A \times \R^{n-1})$ for Borel sets $A \subset \R$. 
For $i = 2,\dots,n$, let $\mu_i(dx_i|x_{1:i-1})$ denote the conditional 
distribution of the $i$-th coordinate under $\mu$ knowing the first $i-1$ 
coordinates $x_1,\dots,x_{i-1}$. Under mild regularity assumptions on $\mu$, 
all these conditional measures are well-defined, and we have a general 
formula for the ``full expectation"
\be
\int h(x)\,d\mu(x) = \int h(x_1,\dots,x_n)\,\mu_n(dx_n|x_{1:n-1}) \dots 
\mu_2(dx_2|x_1) \mu_1(dx_1),
\en
for any bounded measurable function $h$ on $\R^n.$
For example, it suffices to require that $\mu$ has a smooth positive density,
which is polynomially decaying at infinity. Then we will say that $\mu$ is 
regular. In many inequalities, the regularity assumption is only technical 
for purposes of the proof, and may easily be omitted in the resulting 
formulations.

The distance functionals $D$, $I$, and $\mathcal{T}$ satisfy the following 
tensorisation relations with respect to product measures similarly to (5.1). 
To emphasize  the dimension, we denote by $\gamma_n$ the standard Gaussian 
measure on $\R^n$.

\vskip5mm
{\bf Lemma 5.1.} {\it For any regular probability measure $\mu$ on 
$\R^n$ with finite second moment, 
\bee
D(\mu|\gamma_n)
 & = & 
D(\mu_1|\gamma_1) + \sum_{i=2}^n \int D(\mu_i(\,\cdot \ |x_{1:i-1})\, | \gamma_1)\,d\mu(x), \\
I(\mu|\gamma_n)
 & \geq & 
I(\mu_1|\gamma_1) + \sum_{i=2}^n \int I(\mu_i(\,\cdot \ |x_{1:i-1})\, | \gamma_1)\,d\mu(x), \\
\mathcal{T}(\mu,\gamma_n)
 & \leq & 
\mathcal{T}(\mu_1,\gamma_1) + \sum_{i=2}^n 
\int \mathcal{T}(\mu_i(\,\cdot \ |x_{1:i-1}),\gamma_1)\,d\mu(x).
\ene
}

Note that this statement remains to hold also for other product references 
measures $\nu^n$ on $\R^n$ in place of $\gamma_n$ (with necessary regularity
assumptions for the case of Fisher information).

Applying the first two inequalities, we see that the deficit $\delta$ 
satisfies a similar property,
\be
\delta(\mu) \, \geq \,
\delta(\mu_1) + \sum_{i=2}^n \int \delta(\mu_i(\,\cdot \ |x_{1:i-1}))\,d\mu(x).
\en

\vskip5mm
{\bf Proof of Lemma 5.1.} The equality for the relative entropy is 
a straightforward calculation. We refer to Appendix A of [G-L] for 
a (general) tensorisation inequality for transport costs. Below, we sketch 
the proof of the inequality involving Fisher information.

Let $\mu$ be a regular probability measure on $\R^n$ admitting a smooth 
density $f$ with respect to $\gamma_n$. Note that the first marginal 
$\tilde{\mu}$ of $\mu$ on the first $n-1$ coordinates has density $\tilde{f}(x_{1:n-1}) = \int f(x_{1:n-1},x_n)\,\gamma(dx_n)$ and that 
$\mu_n(\,\cdot\,| x_{1:n-1})$ has density 
$f(x_n | x_{1:n-1}) = f(x_{1:n-1},x_{n}) / \tilde{f}(x_{1:n-1})$.
We have
\bee
I(\mu | \gamma_n) 
 & = &
\sum_{i=1}^{n-1} \int \frac{(\partial _{x_i} f)^2}{f} (x)\,\gamma_n(dx) + 
\int \frac{(\partial _{x_n} f)^2}{f} (x)\,\gamma_n(dx)\\
 & = & 
\sum_{i=1}^{n-1} \int \left(\int \frac{(\partial _{x_i} f)^2}{f} (x_{1:n-1},x_n)\,\gamma_1(dx_n)\right)\,\gamma_{n-1}(dx_{1:n-1}) \\
 & & + \ 
\int I(\mu_n(\,\cdot\,| x_{1:n-1})|\gamma_1) \,\tilde{\mu}(dx_{1:n-1})\\
 & \geq & 
\sum_{i=1}^{n-1} \int \frac{(\partial _{x_i} \tilde{f})^2}{\tilde{f}} (x_{1:n-1})\,\gamma_{n-1}(dx_{1:n-1}) + 
\int I(\mu_n(\,\cdot\,| x_{1:n-1})|\gamma_1) \,\tilde{\mu}(dx_{1:n-1})\\
 & = &
I(\tilde{\mu} | \gamma_{n-1}) + 
\int I(\mu_n(\,\cdot\,| x_{1:n-1})|\gamma_1) \,d\mu(x),
\ene
where the inequality holds by an application of Jensen's inequality with 
the function $\psi(u,v) = u^2/v$ which is convex on the upper half-plane 
$\R \times (0,\infty)$. 
The proof is completed by induction. \qed

\vskip5mm
{\bf Proof of Theorem 1.3.}
Let us apply the one dimensional result (4.8) with constant 
$c = 1/(128\,\pi^2)$ in (5.2) to the measures 
$\mu_1$ and $\mu_i(\,\cdot \ |x_{1:i-1})$. Put $t_1 = \int x_1\,\mu_1(dx_1)$, 
$$
t_i(x) = t_i(x_1,\dots,x_{i-1}) = \int x_i\,\mu_i(dx_i | x_{1:i-1}), \qquad
x = (x_1,\dots,x_n) \in \R^n,
$$
and denote by $\tilde \mu_i(\,\cdot \ | x_{1:i-1})$ the corresponding shift 
of $\mu_i(\,\cdot \ | x_{1:i-1})$ as in Corollary 4.3: 
$\tilde \mu_i(\,\cdot \ | x_{1:i-1}) = \mu_i(\,\cdot \ | x_{1:i-1})_{-t_i}$.  
Then we have
$$
256\,\pi^2 \delta(\mu) \, \geq \, 
\frac{\mathcal{T}^2(\tilde \mu_1,\gamma_1)}{D(\tilde \mu_1|\gamma_1)} +
\sum_{i=2}^n \int
\frac{\mathcal{T}^2(\tilde \mu_i(\,\cdot \ | x_{1:i-1}),\gamma_1)}
{D(\tilde \mu_i(\,\cdot \ |x_{1:i-1})|\gamma_1)} \ d\mu(x).
$$
By Jensen's inequality with the convex function $\psi(u,v) = u^2/v$ 
($u \in \R$, $v \geq 0$), 
\bee
256\,\pi^2 \delta(\mu) 
 & \geq &
\frac{\mathcal{T}^2(\tilde \mu_1,\gamma_1)}{D(\tilde \mu_1|\gamma_1)} + 
\sum_{i=2}^n
\frac{\big(\int \mathcal{T}(\tilde \mu_i(\,\cdot \ | x_{1:i-1}),\gamma_1)\, d\mu(x)\big)^2} 
{\int D(\tilde \mu_i(\,\cdot \ |x_{1:i-1})|\gamma_1) \, d\mu(x)} \\
 & \geq &
\frac{\big(\mathcal{T}(\tilde \mu_1,\gamma_1) + \sum_{i=2}^n 
\int \mathcal{T}(\tilde \mu_i(\,\cdot \ | x_{1:i-1}),\gamma_1)\,d\mu(x)\big)^2}
{D(\tilde \mu_1|\gamma_1) + \sum_{i=2}^n 
\int D(\tilde \mu_i(\,\cdot \ |x_{1:i-1})|\gamma_1) \, d\mu(x)},
\ene
where the last bound comes from the inequality
$$
\sum_{i=1}^n \psi(u_i,v_i) \geq 
\psi\Big(\sum_{i=1}^n u_i,\sum_{i=1}^n v_i\Big),
$$
which is due to the convexity of $\psi$ and its 1-homogeneity. Note that 
the first inequality could also be proved by using Cauchy-Schwarz inequality.

Now consider the map $T:\R^n\to\R^n$ defined for all $x\in \R^n$ by 
$$
T(x) = \big(x_1-t_1,x_2-t_2(x_1),\ldots,x_n-t_n(x_1,x_2,\ldots,x_{n-1})\big).
$$
By definition, $T$ pushes forward $\mu$ onto $\bar{\mu}$. The map $T$ is 
invertible and its inverse $U= (u_1,\dots,u_n)$ satisfies
\bee
u_1(x) & = & x_1 + t_1, \\
u_2(x) & = & x_2 + t_2(u_1(x)),  \\
       & \vdots & \\ 
u_i(x) & = & x_i + t_i(u_1(x),\dots, u_{i-1}(x)),  \\   
       & \vdots & \\ 
u_n(x) & = & x_n + t_n(u_1(x),\dots, u_{n-1}(x)).    
\ene
It is not difficult to check that
$\bar{\mu}_1=\tilde{\mu}_1$ and for all $i\geq 2$, $\bar{\mu}_i(\,\cdot\, | x_{1:i-1}) = \tilde{\mu}_i(\,\cdot\, | u_1(x),\ldots,u_{k-1}(x)).$
%Set $\bar \mu_1 = \tilde \mu_1$ and  
%$\bar \mu_i(\,\cdot \ |x_{1:i-1}) = 
%\tilde \mu_i(\,\cdot \ |u_1(x),\dots, u_{i-1}(x))$
%for $i \geq 2$. Finally define the probability measure $\bar \mu$ 
%by virtue of the identity such as (5.1):
%$$
%\int f(x)\,\bar d\mu(x) = \int f(x_1,\dots,x_n)\,\bar \mu_n(dx_n|x_{1:n-1}) \dots 
%\bar \mu_2(dx_2|x_1) \bar \mu_1(dx_1).
%$$
%It then should be clear that the map $U:\R^n \rightarrow \R^n$
%pushes forward $\bar{\mu}$ onto $\mu$. 
Therefore, since $U$ pushes forward $\bar{\mu}$ onto $\mu$, 
$$
\hskip-45mm \mathcal{T}(\tilde \mu_1,\gamma_1) + \sum_{i=2}^n 
\int \mathcal{T}(\tilde \mu_i(\,\cdot \ | x_{1:i-1}),\gamma_1) \, d\mu(x) \\ 
$$
$$
\hskip22mm = \ 
\mathcal{T}(\bar \mu_1,\gamma_1) + \sum_{i=2}^n 
\int \mathcal{T}(\tilde \mu_i(\,\cdot \ |u_1(x),\dots, u_{i-1}(x)),\gamma_1) \,  
d\bar \mu(x)
$$
$$
\hskip25mm = \ 
\mathcal{T}(\bar \mu_1,\gamma_1) + \sum_{i=2}^n 
\int \mathcal{T}(\bar \mu_i(\,\cdot \ | x_{1:i-1}),\gamma_1)\,d\bar\mu(x) \ \geq \
\mathcal{T}(\bar \mu,\gamma_n),
$$
where we made use of Lemma 5.1 on the last step. The same with equality 
sign holds true for the $D$-functional. As a result, in terms of the 
recentered measure $\bar \mu$, we arrive at the following bound:
\be
D(\mu|\gamma_n) + \frac{1}{256\,\pi^2}\,
\frac{\mathcal{T}^2(\bar \mu,\gamma_n)}{D(\bar \mu|\gamma_n)} \leq \frac{1}{2}\,I(\mu|\gamma_n).
\en
%To end the proof of Theorem 1.3, it remains to recognize $\bar \mu$ as the 
%distribution of the recentered random vector $\bar X$ (which follows from the above
%construction). 
Thus, we have established in (5.3) the desired inequality (1.10) 
with constant $c = \frac{1}{256\,\pi^2}$.
\qed

\vskip5mm
{\bf Remark 5.2.}
In order to relate the transport distance $\mathcal{T}$ to $W_1$, one may
apply Lemma 2.3. Following the very definition of the transport distances,
it implies that
$$
(1 - \log 2)\,\min\{W_1(\mu,\nu),W_1^2(\mu,\nu)\} \, \leq \, 
\mathcal{T}(\mu,\nu) \, \leq \, W_1(\mu,\nu),
$$
for all probability measures $\mu$ and $\nu$ on $\R^n$.

\vskip5mm
The proof of Theorem 1.4 will make use of the classical Pr\'ekopa-Leindler 
theorem, which we state below.

\vskip5mm
{\bf Theorem 5.3.} {([Pr1, Pr2], [Le]) 
{\it For a number $t\in (0,1)$, assume that measurable functions 
$f,g,h: \R^d \to \R$ satisfy 
$$
h((1-t)x + t y) \leq (1-t) f(x) + t g(y), \qquad for \ all \ x,y \in \R^d.
$$
Then 
$$
\int e^{-h(z)}\,dz \geq \left(\int e^{-f(x)}\,dx\right)^{1-t} \left(\int e^{-g(y)}\,dy\right)^t.
$$
}

\vskip5mm
{\bf Proof of Theorem 1.4.} It is similar to the proof of Theorem 1.3. 
The main point is that, if $\mu$ has a smooth density $f=e^{-V}$ with respect 
to Lebesgue measure, with a $V$ such that $V''\geq \ep\,{\rm I_n}$ for some 
$\ep>0$, then the first marginal $\mu_1$ has a density of the form $e^{-v_1}$ 
with $v_1''\geq \ep$. Moreover, for each $i = 2,\dots,n$ and all $x\in \R^n$, 
the one dimensional conditional probability $\mu_i(\,\cdot\,|x_{1:i-1})$ 
has a density $e^{-v_i(x_i|x_{1:i-1})}$ with 
$\left(\partial^2/\partial x_i^2\right) v_i(x_i|x_{1:i-1})\geq \ep$. 
Indeed, by definition of conditional probabilities,
$$
v_i(x_i|x_{1:i-1}) = 
-\log\left( \int e^{-V(x_{1:i},y_{i+1:n})}\,dy_{i+1}\cdots dy_n\right) +
w(x_{1:i-1}),
$$
where 
$w(x_{1:i-1}) = 
\log\left( \int e^{-V(x_{1:i-1},y_{i:n})}\,dy_idy_{i+1}\cdots dy_n\right)$ 
does not depend on $x_i$. Since $V'' \geq \ep\,{\rm I_n}$, for any 
$i = 2,\dots,n$ and any $x\in \R^n$, the function 
$$
(y_i,y_{i+1},\ldots,y_n) \mapsto 
V(x_{1:i-1},y_i,\ldots,y_n) - \frac{\ep}{2}\,y_i^2
$$
is convex. Thus defining, for $t\in (0,1)$, $x\in \R^n$ and $a_i,b_i \in \R$,
the functions
\bee
f(y_{i+1},\ldots,y_n) 
 & = & 
V(x_{1:i-1}, a_i,y_{i+1:n}) - \frac{\ep}{2}\, a_i^2,\\
g(y_{i+1},\ldots,y_n) 
 & = & 
V(x_{1:i-1}, b_i,y_{i+1:n}) - \frac{\ep}{2}\, b_i^2,\\
h(y_{i+1},\ldots,y_n) 
 & = & 
V(x_{1:i-1}, (1-t)a_i + tb_i,y_{i+1:n}) - 
\frac{\ep}{2} \left((1-t)a_i + t b_i\right)^2,
\ene
one sees that 
$$ 
h((1-t)y_{i+1:n} + t z_{i+1:n}) \leq (1-t) f(y_{i+1:n}) + t g(z_{i+1:n}),
\qquad {\rm for \ all} \ \ y,z\in \R^n.
$$
Therefore, applying Theorem 5.3 to the triple $(f,g,h)$, one gets easily 
that
$$
v_i((1-t)a_i+tb_i|x_{1:i-1}) \leq 
(1-t) v_i(a_i | x_{1:i-1}) +t v_i(b_i | x_{1:i-1}) - 
\frac{\ep}{2}\, t(1-t)(a_i-b_i)^2.
$$
Since $v_i$ is smooth, this inequality is equivalent to 
$\left(\partial/ \partial x_i\right)^2 v_i(x_i|x_{1:i-1})\geq \ep.$ 
A similar conclusion holds for $v_1$. Therefore, $\mu_1$ and the conditional probabilities $\mu_i(\,\cdot\,|x_{1:i-1})$ verify the assumption of 
Corollary 4.4. Thus, applying the tensorisation formula (5.2), we get
$$
\delta(\mu)\geq c\min\{1,\ep\} \bigg(W_2^2(\tilde{\mu}_1,\gamma_1) + 
\sum_{i=2}^n W_2^2(\tilde{\mu}_i(\,\cdot\,|x_{1:i-1}),\gamma_1)\bigg),
$$
where, as before, $\tilde{\mu}_i(\,\cdot\,|x_{1:i-1})$ is the shift of $\mu_i(\,\cdot\,|x_{1:i-1})$ by its mean. Reasoning as in the proof of 
Theorem 1.3, we see that the quantity inside the brackets is bounded from 
below by $W_2^2(\bar{\mu},\gamma_n)$.
\qed

% -----------------------------  section 6 ---------------------------

\vskip10mm
\section{Appendix A: The reversed transport-entropy inequality}
\setcounter{equation}{0}

Here we include a simple proof of the general inequality of Lemma 3.2,
$$
D(X_t |Y_t)\leq \frac{1}{2t}\, W_2^2(X,Y), \qquad t > 0,
$$
where $X$ and $Y$ are  random vectors in $\R^n$ with finite second moments.

We denote by $p_U$ the density of a random vector $U$ and by $p_{U|V=v}$ 
the conditional density of $U$ knowing the value of a random vector $V=v$. 
Note that the regularized random vectors $X_t = X + \sqrt{t}\,Z$ have smooth 
densities.

By the chain rule formula for the relative entropy, one has 
$$
D(X,Y,X_t|X,Y,Y_t) = 
D(X_t|Y_t)+\int D(p_{X,Y|X_t=v}| p_{X,Y|Y_t=v})\, p_{X_t}(v)\, dv,
$$
and therefore
$$
D(X,Y,X_t|X,Y,Y_t)\geq  D(X_t|Y_t).
$$
On the other hand, we also have
$$
D(X,Y,X_t|X,Y,Y_t) = \int\!\!\!\int 
D(p_{X_t|(X,Y)=(x,y)}|p_{Y_t|(X,Y)=(x,y)})\, p_{X,Y}(x,y)\, dx dy.
$$

Now observe that $p_{X_t|(X,Y)=(x,y)}$ is the density of a normal law with mean 
$x$ and covariance matrix $t I_n$, and similarly for $p_{Y_t|(X,Y)=(x,y)}$. 
But
$$
D(x + \sqrt{t}Z\,|\,y + \sqrt{t}\,Z) = \frac{|x-y|^2}{2t},
$$
so 
$$
D(X,Y,X_t|X,Y,Y_t) = \frac{1}{2t} \int\!\!\!\int  |x-y|^2  p_{X,Y}(x,y)\,dx dy
 = \frac{1}{2t}\, W_2^2(X,Y), 
$$ 
where the last equality follows by an optimal choice for the coupling density 
of $X$ and $Y$.

\vskip10mm
\section{Appendix B: Reinforced transport-entropy inequalities}
\setcounter{equation}{0}

In this section, we explain how to derive Theorem 4.1 in the form (4.3).

\vskip5mm
{\bf Proof of Theorem 4.1.} To derive the inequality (4.3) for probability 
measures with mean zero, we follow an argument of [B-K]. Let $\mu$ be 
a probability measure on $\R$ such that $D(\mu |\gamma)$ is finite and 
consider the monotone rearrangement map $T$ transporting $\gamma$ 
onto $\mu$. It is defined by $T(x) = F_{\mu}^{-1} \circ F_{\gamma}(x)$, 
where $F_\mu(x) = \mu (-\infty,x]$ and $F_{\gamma}(x)= \gamma (-\infty,x]$
are the corresponding distribution functions, and
$F_\mu^{-1}(t) = \inf\{x \in \R: F_\mu(x) \geq t\}$
is the generalized inverse of $F_\mu$ (defined for $0 < t < 1$). 
It is well known that $T$ pushes forward $\gamma$ on $\mu$ and achieves 
the minimal value in the optimal transport problem:
$$
W_2^2(\mu,\gamma) = \int (T(x)-x)^2\,d\gamma(x).
$$

The starting point is the following inequality going back to Talagrand's 
paper [T] (see equation (2.5) of [T]): 
\begin{eqnarray}
D(\mu|\gamma)
 & \geq &
\frac{1}{2}\,W_2^2(\mu,\gamma) + 
\int \big(T'(x) - 1 - \log T'(x))\big)\,d\gamma(x) \nonumber \\
 & \geq &
\frac{1}{2}\,W_2^2(\mu,\gamma) + \int \Delta(|T'(x) - 1|\big)\,d\gamma(x),
\end{eqnarray}
where the second inequality comes from the fact that 
$\Delta(x) \geq \Delta(|x|)$ for all $x>-1$.
%Since $V'' \geq \kappa$, the measure $\nu$ may be obtained as the image of 
%the standard Gaussian measure $\gamma$ under an increasing map whose Lipschitz 
%norm $\leq 1/\sqrt{\kappa}$. As a one dimensional statement, it is rather simple 
%and can be verified directly. 
On the other hand, $\gamma$ is known to satisfy 
the Cheeger-type analytic inequality
\be
\lambda \int |f - m(f)|\,d\gamma \leq \int |f'|\,d\gamma
\en
with optimal constant $\lambda = \sqrt{\frac{2}{\pi}}$ 
(see e.g Theorem 1.3 of [B-H]). Here, $f:\R \rightarrow \R$ 
may be an arbitrary locally Lipschitz function with Radon-Nikodym derivative 
$f'$, and $m(f)$ denotes a median of $f$ under $\gamma$. 
%Hence,
%\be
%\lambda \sqrt{\kappa} \int |f - m(f)|\,d\nu \leq \int |f'|\,d\nu
%\en
%with the median functional understood with respect to $\nu$.
According to Theorem 3.1 of [B-H], (7.2) can be generalized as
\be
\int L(f - m(f))\,d\gamma \leq \int L(c_L f'/\lambda)\,d\gamma
\en
with an arbitrary even convex function $L:\R \rightarrow [0,\infty)$, 
such that $L(0)=0$, $L(t)>0$ for $t>0$, and 
$$
c_L = \sup_{t>0} \frac{tL'(t)}{L(t)} < \infty,
$$
where $L'(t)$ may be understood as the right derivative at $t$.

We apply (7.3) with $L(t) = \Delta(|t|) = |t| - \log(1+|t|)$ in which case 
$c_L = 2$, so that
\be
\int \Delta(|f - m(f)|)\,d\gamma \leq \int \Delta(2\, |f'|/\lambda)\,d\gamma.
\en
It will be convenient to replace here the median with the mean
$\gamma(f) = \int f\,d\gamma$. First observe that, by Jensen's inequality, 
(7.4) yields
\be
\Delta(|\gamma(f) - m(f)|) \leq \int \Delta(2\, |f'|/\lambda)\,d\gamma.
\en
Hence, using once more the convexity of $\Delta$ together with (7.4)-(7.5) 
for the function $2f$, we get
\bee
\int \Delta(|f-\gamma(f)|)\,d\gamma
 & \leq &
\frac{1}{2}\int \Delta\big(2\,|f-m(f)|\big)\,d\gamma +
\frac{1}{2}\,\Delta\big(2\,|\gamma(f) - m(f)|\big) \\
 & \leq &
\int \Delta(4\, |f'|/\lambda)\,d\gamma. 
\ene
Equivalently,
$$
\int \Delta(|f'|)\,d\gamma \geq 
\int \Delta\Big(\frac{\lambda}{4}\,|f-\gamma(f)|\Big)\,d\gamma. 
$$
To further simplify, one may use the lower bound $a)$ of Lemma 2.3 which
yields
$$
\int \Delta(|f'|)\,d\gamma \geq \Big(\frac{\lambda}{4}\Big)^2
\int \Delta(|f-\gamma(f)|)\,d\gamma. 
$$
It remains to apply the latter with $f(x) = T(x) - x$ when estimating the last 
integral in (7.1). Since $\mu$ and $\gamma$ have mean zero, this gives
$$
D(\mu|\gamma) \geq \frac{1}{2}\,W_2^2(\mu,\gamma) + \frac{1}{8\pi}
\int \Delta(T(x) - x\big)\,d\gamma(x),
$$
and the last integral is certainly greater than (and actually equals to) $\mathcal{T}(\mu,\gamma)$.
\qed

%\vskip5mm
%{\bf Theorem 8.1.} {\it Let $\nu$ be a probability measure on the real line 
%with density $\frac{d\nu(x)}{dx} = e^{-V(x)}$ such that $V''(x) \geq \kappa$, 
%for all $x \in \R$ with some $\kappa > 0$. Then, for any probability measure 
%$\mu$ on $\R$ with mean zero,
%$$
%D(\mu|\nu) \geq \frac{1}{2}\,W_2^2(\mu,\nu) + 
%\min\big(\lambda \sqrt{\kappa}, (\lambda \sqrt{\kappa})^2\big)\, \mathcal{T}(\mu,\nu),
%$$
%where $\lambda = \sqrt{\frac{1}{8 \pi}}$.
%}
%
%\vskip5mm
%In the Gaussian case $\nu = \gamma$, we have $\kappa = 1$, and the above inequality 
%turns into (4.3).

\vskip5mm
{\bf Proof of Theorem 4.2.} Let us return to 
the inequality (7.1), i.e.,
\be
D(\mu|\gamma) \geq 
\frac{1}{2}\,W_2^2(\mu,\gamma) + \int \Delta(T'(x) - 1\big)\,d\gamma(x).
\en
The basic assumption (4.4) ensures that $T$ has a Lipschitz norm 
$\leq \frac{1}{\sqrt{\ep}}$, so $T'(x) \leq \frac{1}{\sqrt{\ep}}$. 
Using in (7.6) the lower quadratic bounds on $\Delta$ given in 
$b)$ and $c)$ of Lemma 2.3, we obtain that
\be
D(\mu|\gamma) \geq 
\frac{1}{2}\,W_2^2(\mu,\gamma) + c(\ep) \int (T'(x) - 1)^2\,d\gamma(x),
\en
where
$$
c(\ep) = \frac{1}{2}, \ \ {\rm for} \ \ep \geq 1, \qquad \quad
c(\ep) = \frac{\Delta(\frac{1}{\sqrt{\ep}} - 1)}{(\frac{1}{\sqrt{\ep}} - 1)^2}, \ \ 
{\rm for} \ 0 < \ep < 1.
$$

On the other hand, applying the Poincar\'e-type inequality for the Gaussian 
measure
$$
\Var_\gamma(f) \leq \int f'^2\,d\gamma
$$
with $f(x) = T(x) - x$, together with the assumption that
$\int x\,d\mu(x) = \int T(x)\,d\gamma(x) = 0$, the last integral in (7.7)
can be bounded from below by
$$
\int (T(x) - x)^2\,d\gamma(x) = W_2^2(\mu,\gamma).
$$
It remains to use, for $0 < \ep < 1$, the bound 
$\Delta(a) \geq (1 - \log 2)\,\min\{a,a^2\}$.
The inequality (4.5) is proved.
\qed

\vskip10mm
\section{Appendix C: Equality cases in the logarithmic Sobolev inequality 
for the standard Gaussian measure}
\setcounter{equation}{0}

In this last section, we show how Theorem 1.3 can be used to recover 
the following result by E. Carlen [C].

\vskip5mm
{\bf Theorem 8.1.} {([C]) \it Let $\mu$ be a probability measure on $\R^n$ 
such that $D(\mu|\gamma) < \infty$. We have
$$
D(\mu|\gamma) = \frac{1}{2}\,I(\mu|\gamma),
$$ 
if and only if $\mu$ is a translation of $\gamma.$
}

\vskip5mm
In what follows, we denote by $\mathcal{S}_n$ the set of permutations of $\{1,\ldots,n\}$. If $\mu$ is a probability measure on $\R^n$, we denote by $\mu_\sigma$ its image under the permutation map
$$
(x_1,\ldots,x_n)\mapsto (x_{\sigma(1)},\ldots,x_{\sigma(n)}).
$$
If $\mu$ has density $f$ with respect to the standard $n$-dimensional 
Gaussian measure $\gamma$, then the density of $\mu_\sigma$ with respect 
to $\gamma$ is given by
$$
f_\sigma(x_1,\ldots,x_n)=f(x_{\sigma^{-1}(1)},\ldots,x_{\sigma^{-1}(n)}).
$$
Obviously,
$$
I(\mu_\sigma|\gamma) = I(\mu|\gamma) \quad \text{and} \quad 
D(\mu_\sigma|\gamma) = D(\mu|\gamma).
$$
Hence, we have the following automatic improvement of Theorem 1.3.

\vskip5mm
{\bf Theorem 8.2.} {\it Let $X$ be a random vector in $\R^n$ with law $\mu$.
Then,
$$
D(\mu |\gamma) + c\max_{\sigma \in \mathcal{S}_n}
\frac{\mathcal{T}^2(\overline{\mu_\sigma},\gamma)}{D(\overline{\mu_\sigma}|\gamma)}
\leq \frac{1}{2}I(\mu |\gamma),
$$
where $\overline{\mu_\sigma}$ is the law of the random vector $Y^{\sigma}$ 
defined by
$$
Y^\sigma_i =
X_{\sigma(i)} - \E(X_{\sigma(i)}| X_{\sigma(1)},\ldots,X_{\sigma(i-1)}).
$$
}

%\vskip5mm
%Thanks to this result, we can recover Carlen's characterization of equality cases in %the logarithmic Sobolev inequality for the standard Gaussian on $\R^n.$ 
\vskip5mm
{\bf Proof of Theorem 8.1.}
To avoid complicated notations, we will restrict ourselves to the dimension 
$n=2$. We may assume that $\mu$ has a smooth density $p$ with respect to the
Lebesgue measure such that $D(\mu |\gamma) = \frac{1}{2}\,I(\mu |\gamma) < \infty$. 
Necessarily, $\mu$ has a finite second moment, and moreover, 
$\overline{\mu_\sigma} = \gamma$, for all $\sigma \in \mathcal{S}_2$, i.e., 
for $\sigma = \mathrm{id} = (12)$ and $\sigma = (21)$.

For a random vector $X$ with law $\mu$, put $m_1 = \E X_1$, $m_2 = \E X_2$, 
$a(X_1) = \E\,(X_2|X_1)$ and $b(X_2) = \E\,(X_1|X_2)$.
The probability measure $\gamma = \overline{\mu_\mathrm{id}}$ represents 
the image of $\mu$ under the map $(x_1,x_2) \mapsto (x_1-m_1,x_2-a(x_1))$. 
It then easily follows that
$$
p(x_1,x_2) = \frac{1}{2\pi}\exp\left(-\frac{1}{2}\, (x_1-m_1)^2 - \frac{1}{2}\,(x_2-a(x_1))^2\right)
$$
for almost all $(x_1,x_2) \in \R^2$. Since also 
$\gamma = \overline{\mu_{(2,1)}}$, the same reasoning yields
$$
p(x_1,x_2) = \frac{1}{2\pi}\exp\left(-\frac{1}{2}\, (x_2-m_2)^2 - \frac{1}{2}\,(x_1-b(x_2))^2\right),
$$
for almost all $(x_1,x_2) \in \R^2$. Therefore, for almost all 
$(x_1,x_2)\in \R^2$, it holds
$$
(x_1-m_1)^2 +(x_2-a(x_1))^2=(x_2-m_2)^2 +(x_1-b(x_2))^2.
$$
Let us denote by $A$ the set of all couples $(x_1,x_2)$ for which there 
is equality, and for $x_1\in \R$, let $A_{x_1} = \{x_2 \in \R: (x_1,x_2)\in A\}$
denote the corresponding section of $A$. By Fubini's theorem, 
$$
0 = |\R^2 \setminus A| = \int_{-\infty}^\infty |\R \setminus A_{x_1}|\,dx_1,
$$
where $| \cdot |$ stands for the Lebesgue measure of a set in the corresponding 
dimension. Hence, for almost all $x_1$, the set $\R\setminus A_{x_1}$ is 
of Lebesgue measure $0$. For any such $x_1$,
$$
2x_2(m_2-a(x_1)) + a(x_1)^2 -m_2^2 + (x_1-m_1)^2\geq 0,\qquad 
\forall \, x_2\in A_{x_1}.
$$
Thus, $a(x_1)=m_2$ (otherwise letting $x_2\to\pm \infty$ would lead to 
a contradiction). This proves that $a=m_2$ almost everywhere, and therefore, 
the random vector $(X_1-\E X_1, X_2-\E X_2)$ is standard Gaussian. But this 
means that $\mu$ is a translation of $\gamma$. 
\qed

\vskip5mm
{\bf Acknowledgement.} We would like to thank M. Ledoux for interesting 
comments and for pointing to the paper by F-Y. Wang. We also thank an anonymous referee for valuable suggestions and for pointing out a mistake in the initial version of this paper.

\vskip1cm


\begin{thebibliography}{}

\end{thebibliography}


\vskip5mm
\begin{thebibliography}{BH3}
\small

\bibitem[A]{A} An{\'e}, C., Blach{\`e}re, S., Chafa{\"{\i}}, D., Foug{\`e}res, P., Gentil, I.,
  Malrieu, F., Roberto, C., Scheffer, G.
Sur les in\'egalit\'es de {S}obolev logarithmiques. volume~10
  Panoramas et Synth\`eses [Panoramas and Syntheses].
Soci\'et\'e Math\'ematique de France, Paris, 2000. With a preface by Dominique Bakry and Michel Ledoux.

\bibitem[B-B-G]{B-B-G}
Bakry, D., Bolley, F., Gentil, I. Dimension dependent hypercontractivity 
         for Gaussian kernels. Probab. Theory Related Fields 154 (2012), 
         no. 3-4, 845--874. 

\bibitem[B-E]{B-E}
Bakry, D., \'Emery, M. Diffusions hypercontractives. Seminaire de probabilites, 
         XIX, 1983/84, 177--206, Lecture Notes in Math., 1123, Springer, 
         Berlin, 1985.

\bibitem[B-L]{B-L} 
Bakry, D., Ledoux, M. A logarithmic Sobolev form of the Li-Yau parabolic 
        inequality. Rev. Mat. Iberoam. 22 (2006), no. 2, 683¬ñ-702. 

\bibitem[B-K]{B-K} 
Barthe, F., Kolesnikov, A. V. Mass transport and variants of the logarithmic 
         Sobolev inequality. J. Geom. Anal. 18 (2008), no. 4, 921ñ-979. 

\bibitem[Bl]{Bl}
Blachman, N. M. The convolution inequality for entropy powers. 
         IEEE Trans. Inform. Theory 11 (1965), 267--271.

\bibitem[B]{B} 
Bobkov, S. G. An isoperimetric inequality on the discrete cube, and an elementary
         proof of the isoperimetric inequality in Gauss space.
         Ann. Probab. 25 (1997), no. 1, 206--214.
         
\bibitem[B-G-L]{B-G-L} 
Bobkov, S. G., Gentil, I., Ledoux, M. Hypercontractivity of Hamilton-Jacobi 
         equations. J. Math. Pures Appl. (9) 80 (2001), no. 7, 669--696.

         
\bibitem[B-G]{B-G} 
Bobkov, S. G., Gotze, F. Exponential integrability and transportation 
         cost related to logarithmic Sobolev inequalities. 
         J. Funct. Anal. 163 (1999), no. 1, 1--28.

\bibitem[B-H]{B-H} 
Bobkov, S. G., Houdr\'e, C. Isoperimetric constants for product probability 
         measures. Ann. Probab. 25 (1997), no. 1, 184--205.

\bibitem[C]{C}     
Carlen, E. A. Superadditivity of Fisher's information and logarithmic 
         Sobolev inequalities. J. Funct. Anal. 101 (1991), no. 1, 194--211.

\bibitem[C-F]{C-F}     
Carlen, E., Figalli, A.
Stability for a GNS inequality and the Log-HLS inequality, with application to the critical mass Keller-Segel equation. Duke Math. J., 162 (2013), no. 3, 579--625.


\bibitem[C-F-M-P]{C-F-M-P}
Cianchi, A., Fusco, N., Maggi, F., Pratelli, A. On the isoperimetric deficit 
         in Gauss space. Amer. J. Math. 133(1):131-186, 2011.

\bibitem[CE]{CE}
Cordero-Erausquin, D. Some applications of mass transport to Gaussian-type 
         inequalities. Arch. Ration. Mech. Anal. 161 (2002), no. 3, 257--269. 

           
\bibitem[D-C-T]{D-C-T} 
Dembo, A., Cover, T. M., Thomas, J. A. Information-theoretic 
         inequalities. IEEE Trans. Inform. Theory 
         37 (1991), no. 6, 1501--1518.

\bibitem[D-T]{D-T}
Dolbeault, J., Toscani, G. Improved interpolation inequalities, relative entropy and fast
  diffusion equations. To appear in Annales de l'Institut Henri Poincare (C) Non Linear
  Analysis, (2013).

\bibitem[E]{E}
Eldan, R. A two-sided estimate for the Gaussian noise stability deficit.
         Preprint (2013), arXiv:1307.2781 [math.PR]. 
         
\bibitem [E-K-S]{E-K-S} 
Erbar, M., Kuwada, K., Sturm, K-T. On the equivalence of the entropic 
         curvature-dimension condition and Bochner's inequality on metric 
         measure spaces. Preprint (2013), arXiv:1303.4382 [math.DG].

\bibitem[F-M-P1]{F-M-P1}
Figalli, A., Maggi, F., Pratelli, A. A refined {B}runn-{M}inkowski inequality for convex sets.
Ann. Inst. H. Poincar\'e Anal. Non Lin\'eaire, 26(6) (2009), 2511--2519.


\bibitem[F-M-P2]{F-M-P2}
Fusco, N., Maggi, F., Pratelli, A. The sharp quantitative Sobolev inequality for functions
of bounded variation. J. Func. Anal., 244, (2007), 315–341.

\bibitem[G-L]{G-L}
Gozlan, N., L\'eonard, C. Transport inequalities - A survey. Markov Processes 
         and Related Fields 16 (2010), 635--736.

\bibitem[G]{G} 
Gross, L. Logarithmic Sobolev inequalities. Amer. J. Math. 97 (1975), 
         1061--1083.

\bibitem[I-M]{I-M}
Indrei, E., Marcon, D. A quantitative log-Sobolev inequality for a two parameter
         family of functions. To appear in Int. Math. Res. Not. (2013).

\bibitem[L1]{L1}  
Ledoux, M. Concentration of measure and logarithmic Sobolev inequalities.
         Seminaire de Probabilites XXXIII. Lecture Notes in Math. 1709 (1999), 
         120-216, Springer. 
         
\bibitem[L2]{L2}           
Ledoux, M. The concentration of measure phenomenon. Math. Surveys and monographs,
         vol. 89, AMS, 2001.
         
\bibitem[Le]{Le} 
Leindler, L., On a certain converse of H\"older's inequality II, stochastic
programming, Acta Sci. Math. Szeged 33 (1972), 217--223.

\bibitem[Li]{Li}
Lieb, E. H. Proof of an entropy conjecture of Wehrl. 
         Comm. Math. Phys. 62 (1978), no. 1, 35--41.

%\bibitem[MC]{MC}
%McCann, R. J. Existence and uniqueness of monotone measure-preserving maps. 
%         Duke Math. J. 80 (1995), no. 2, 309¬ñ-323. 

\bibitem[M-N]{M-N}
Mossel, E., Neeman, J. Robust dimension free isoperimetry in Gaussian space.
         Preprint (2012). To appear in Ann. Probab.
         
\bibitem[P]{P}      
Pinsker, M. S. Information and information stability of random variables 
         and processes. Translated and edited by Amiel Feinstein Holden-Day, 
         Inc., San Francisco, Calif.-London-Amsterdam, 1964, xii+243 pp.
         
\bibitem[Pr1]{Pr1}
Pr\'ekopa, A., Logarithmic concave measures with applications to stochastic
         programming. Acta Sci. Math. Szeged 32 (1971), 301--316.

\bibitem[Pr2]{Pr2}
Pr\'ekopa, A., On logarithmic concave measures and functions. 
         Acta Sci. Math. Szeged 34 (1973), 335--343.

\bibitem[O-V]{O-V}
Otto, F., Villani, C. Generalization of an inequality by Talagrand, and links 
         with the logarithmic Sobolev inequality. J. Funct. Anal. 173 (2000), 
         361--400.

\bibitem[R-S]{R-S}  
Raginsky, M., Sason, I. Concentration of measure inequalities in 
         Information Theory. Communications and Coding. Foundations and 
         Trends in Communications and Information Theory, vol. 10 (2013), 
         issues 1 and 2, 1--246.

\bibitem[Se]{Se}
Segal, A. Remark on stability of Brunn–Minkowski and isoperimetric inequalities for convex bodies.
Geometric Aspects of Functional Analysis, Lecture Notes in Mathematics, Volume 2050, (2012), 381--391.

\bibitem[St]{St}
Stam, A. J. Some inequalities satisfied by the quantities of information 
         of Fisher and Shannon. Information and Control 2 (1959), 101--112.

\bibitem[T]{T}
Talagrand, M. Transportation cost for Gaussian and other product measures.
         Geom. Funct. Anal. 6 (1996), 587--600.
         
\bibitem[V]{V}
Villani, C. Optimal transport: Old and new. Volume 338 of Grundlehren der 
         Mathematischen Wissenschaften [Fundamental Principles of Mathematical 
         Sciences]. Springer-Verlag, Berlin, 2009. 
         
\bibitem[W]{W}
Wang, F-Y. Generalized transportation-cost inequalities and applications.
         Potential Anal. 28 (2008), no. 4, 321--334.
         
\bibitem[Wu]{Wu}
Wu, Y. A Simple Transportation-Information Inequality, with Applications to HWI 
         Inequalities and Predictive Density Estimation. Preprint 2011.

\end{thebibliography}
\end{document}